\documentclass{interact}
\usepackage{fancyhdr}
\usepackage{amsthm}
\usepackage{etoolbox}
\usepackage{verbatim}
\usepackage{enumerate}
\usepackage{amsmath}
\usepackage{algorithmicx}
\usepackage{algorithm}
\usepackage{algpseudocode}
\usepackage{amssymb}
\usepackage{tikz}

\begin{document}

\title{Summations of Linear Recurrent Sequences}

\author{
\name{Andrew Lohr\textsuperscript{a}\thanks{CONTACT Andrew Lohr. Email: ajl213@math.rutgers.edu}}
\affil{\textsuperscript{a}Department of Mathematics, Rutgers University, 110 Frelinghuysen
Rd., Piscataway, NJ 08854-8019, USA.}
}

\maketitle

\begin{abstract}
We give an extension of Sister Celine's method of proving hypergeometric sum identities that allows it to handle a larger variety of input summands. In particular, we extend the summand to powers of a C-finite sequence times a hypergeometric term. We then apply this to several problems. Some of these applications give new results, and some reprove already known results in an automated way.
\end{abstract}


{\bf Maple packages and data files}
\\\\
This article is accompanied by the Maple package {\tt RecSum.txt}
which is available with an appendix and several example inputs and outputs at

{\tt http://www.math.rutgers.edu/\~{}ajl213/DrZ/Celine/readme.html} 

{\bf Background}
\\\\
One of the earliest steps in automatically proving identities dates back to Sister Mary Celine Fasenmyer's 1945 Ph.D. thesis \cite{F1}. She gave a technique for computing sums of hypergeometric terms, also see \cite{F2}. Her technique concerns sequences of the form  $x_n = \sum_k H(n,k)$, where the sum is over all $k$ so that $H(n,k)$ is non-zero. Because it is summing over all of these $k$, the problems that it can be applied to only make sense if for each $n$ there are only finitely many values of $k$ that cause $H(n,k)$ to be non-zero. Many expressions constructed from binomial coefficients fit these requirements. It also requires that $H(n,k)$ is doubly hypergeometric, meaning both $\frac{H(n,k+1)}{H(n,k)}$ and $\frac{H(n+1,k)}{H(n,k)}$ are rational functions in $n$ and $k$. In order to determine if there is an order $I$ recurrence for the sequence, her technique picks some $J$ and considers 

\[
0 = \sum_{i=0}^I \sum_{j=0}^J y_{i,j}(n) H(n+i,k+j),
\]
where $y_{i,j}(n)$ is an as yet unknown rational function of $n$. If the value picked for $J$ was not large enough then this procedure will fail, and a higher value for $J$ would be considered. Then, by $H$ being hypergeometric, it is able to reduce all of the $H(n+i,k+j) = G_{i,j}(n,k)H(n,k)$ where $G_{i,j}$ is some rational function of $n$ and $k$. From there, divide everything through by $H(n,k)$. Now, we have something of the form 

\[
0 = \sum_{i=0}^I \sum_{j=0}^J G_{i,j} (n,k) y_{i,j}(n).
\]

Combining denominators on the right hand side, and multiplying through by the common denominator, we get that the right hand side becomes a polynomial in $n$ and $k$, with $\{y_{i,j}(n)\}$ thrown in as well. Collect terms by what power of $k$ appears, and then solve for what the $\{y_{i,j}(n)\}$ have to be in order to make all of the coefficients of powers of $k$ equal to zero. We may get unlucky and have no solution, then, we would need to try a larger $I$ to begin with. If however, we find a solution, we plug that into where we first introduced $y_{i,j}(n)$. Since these have no $k$'s in them, and $x_n$ is obtained by summing over all values of $k$ that make the summand nonzero, we have

\[
0 = \sum_{i=0}^I \sum_{j=0}^J y_{i,j}(n) H(n+i,k+j) = \sum_{i=0}^I \left( \sum_{j=0}^J y_{i,j}(n)\right) x_{n+i}.
\]

Which we may write in shift operator notation as 
\[
0 = \left(\sum_i^I \left( \sum_{j=0}^J y_{i,j}(n)\right) N^i\right)x_n.
\]

At this point we say that we are done. First, having a recurrence allows us to compute the sequence out to very large values very quickly, storing only a constant number of terms. Also, once we have a rational recurrence like this for $x_n$ then we can extract as good asymptotics as desired like using techniques by Birkhoff-Trjizinski which has been nicely summarized in \cite{WZ}. Sometimes we are able to recover a really nice formula. Some of the recurrences found by our procedure are very complicated, so there is little hope to always be able to recover a formula.

For a more complete explanation of Sister Celine's method, look at Chapter 4 of \cite{PWZ}. There are some generalizations of Sister Celine's method given in \cite{Z1}, in particular to certain classes of multiple summations and to a continuous analog.

Some of our applications of the expanded method presented in this paper relate to binomial transforms of functions. There are nice treatments of binomial transforms of Fibonacci like sequences given in \cite{Sp}.

{\bf Main result}
\\\\

We take the described technique of Sister Celine and extend it to allow many more kinds of summands. In particular, the sequence can be of the form $x_n =\sum_k^n a_k^d H(k,n)$ where d is any number,  H is hypergeometric, and $a_k$ is some sequence defined by a rational recurrence relation. Since so many sequences can be so described by rational recurrence relations, this is a significant extension in scope. 

It works very similarly to Sister Celine's method, in that we will consider ratios of successive terms. That is,  to find a recurrence with order at most $I$, start with

\[
\sum_{i=0}^I \sum_{j=0}^{J} \frac{H(n+i,k+j)}{H(n,k)} a_{j+k}^d y_{i,j}(n).
\]

Let $D$ be the order of the recurrence describing $\{a_k\}$. Then, we use that relation to rewrite all of the $\{a_{k+j}\}_{j=D}^J$ in terms of $\{a_{k+j}\}_{j=0}^{D-1}$. That is, by repeatedly applying the relation, we can write each $a_{j+k}$ as a linear combination:

\[
a_{j+k}  = \sum_{m=0}^{D-1} c_{k,j,m} a_{k+m},
\]
where for the $j<D$, we just let $c_{k,j,m} = \begin {cases} 1 & j=m\\ 0 & j \neq m \end{cases}$. Then, since we have an expression with $D$ terms to the $d$, we can expand that out to get at most $Dd$ terms. Then, unlike in Sister Celine's method, where we have a polynomial in $k$, we now have a polynomial in $\{k,a_k,a_{k+1},\ldots a_{k+D-1}\}$. But, once we have collected the coefficients of each of the combinations of those variables, we set all of them equal to zero, and then try to solve for the $y_{i,j}(n)$. As in Sister Celine's method, we are not guaranteed that we can find such a solution for our particular choice of $I$ and $J$. We are guaranteed by WZ theory \cite{PWZ} that for a large enough choice of $I$ and $J$, it gives us a recurrence relation that looks like

\[
 0 = \left( \sum_{i=0}^I \left(\sum_{j=0}^J y_{i,j}(n)\right) N^i \right) x_n.
\]

Our whole technique is implemented in a maple package whose address is given at the beginning of this paper. The usefulness of our technique comes from being easily carried out by a computer, since the systems of equations involved quickly get too large for a person. We invite the reader to use this package the next time that the come across a type of summation problem that they want to analyze.

{\bf Application to enumerating chess king walks}
\\\\

Suppose that there is a king wandering around on an infinite $d$-dimensional chess board. We want to know how many of the $(3^d-1)^n$ walks of length $n$ that the king could take would end up bringing him back to where he started. Given a polynomial $p$, we use the notation $Ct(p)$ to denote the constant term of $p$. Then, by using the powers of $z_i$ to keep track of our total displacement in the $i$ dimension, we have:
\begin{align*}
x_n &=Ct\left(\left(\left(\prod_{i=1}^d z_i +z_i^{-1} +1\right) -1\right)^n\right)\\
&=Ct\left(\sum_{k=0}^n\left(\prod_{i=1}^d z_i +z_i^{-1} +1\right)^k \binom{n}{k} (-1)^{n-k}\right)\\
&=\sum_{k=0}^n Ct\left(\left(\prod_{i=1}^d z_i +z_i^{-1} +1\right)^k\right) \binom{n}{k} (-1)^{n-k}\\
&=\sum_{k=0}^n Ct\left(\left(z +z^{-1} +1\right)^k\right)^d \binom{n}{k} (-1)^{n-k}.
\end{align*}
 
Luckily for us, $Ct\left(\left(z +z^{-1} +1\right)^k\right)$ is already well understood. It is the sequence of central trinomial coefficients (A002426 \cite{Sl}). Also luckily, it is known that this sequence satisfies the recurrence
\[
 0 = \left(N^2 - \frac{2n-1}{n}N - \frac{3n-3}{n}\right)x_n.
\]

So, we are in exactly the set up of our technique. In which case, if we let $a_k = Ct\left(\left(z +z^{-1} +1\right)^k\right)$, we can describe the number of $d$ dimensional king walks which end at the origin after taking $n$ steps by 

\[\sum_{k=0}^n a_k^d \binom{n}{k}(-1)^{n-k}.\]

Once the counting problem has been rewritten as this sum, it clearly falls into the scope of our technique. Using it we are able to find rational recurrences (effectively solve) for all dimensions up to 4.
For a two dimensional king walking around, if we let
\begin{align*}
g(n,N) =& (3n^3+40n^2+175n+250)N^3  \\
 &+(9n^3+138n^2+703n+1190)N^2\\
&+(108n^3+1548n^2+7364n+11632)N\\
&+96n^3+1280n^2+5632n+8192
\end{align*}
then
\[
0 = g(n,N)x_n.
\]
Although this recurrence already looks a little ugly, at least it is short, which is more than can be said of those expressions describing higher dimensions. But they are included in an appendix. Also important is that they were found by a computer.

Something probably more insightful than these walls of text that exactly describe these sequences is their asymptotics:

For the two dimensional king, the number of paths of length n is:
\[
c_2 \frac{8^n}{n}\left( 1 - \frac{4}{9n} + \frac{1}{18n^2} + O\left(\frac{1}{n^3}\right)\right),
\]
For three dimensions the number is:
\[
c_3\frac{26^n}{n^\frac{3}{2}}\left( 1 - \frac{11}{18n} + \frac{683}{5832n^2} + O\left(\frac{1}{n^3}\right)\right).
\]
and for four dimensions the number is:
\[
c_4 \frac{80^n}{n^2} \left(1 - \frac{25}{9n} + \frac{36439}{6561n^2}+O\left(\frac{1}{n^3}\right)\right).
\]

The dominant asymptotics are somewhat unsurprising. The exponential part is all possible paths. The dominant power of $n$ is $\left(\frac{1}{\sqrt{n}}\right)^d$. It is well known that the central binomial coefficient is asymptotically $\frac{2^n}{\sqrt{n}}$, and we are doing something somewhat like that in $d$ dimensions. The value of $c_2$ is approximately equal to $\frac{2}{3\pi}$. This value for $c_2$ can be proven in a rigorous way using classical analysis. For $c_3$ and $c_4$, we are not so lucky, instead, all we can say from non-rigorous observation is that $c_3 \approx .110225343716$ and $ c_4 \approx .068412392872$. There might be some way using a more traditional approach that would get us the true value of these constants.

The $d=2$ case was first worked out by a computer using a different approach. For more information on this, see {E}. Their approach expresses the quantity as a double contour integral and applies their own automated techniques to evaluate it. For information on the techniques, see \cite{AZ}. A completely human produced analysis of this sequence proves more illusive.

{\bf Application to other sequences}
\\\\

Our technique also allows for computing binomial transforms of interesting sequences. An example of this is if we were to let $F_k$ denote the $k$-th Fibonacci number and consider the sequence

\[
x_n = \sum_{k=0}^n F_k \binom{n}{k},
\]
we immediately receive that the recurrence that defines $x_n$ is $0=(-N^2 + 3N - 1)x_n$. This recurrence is identical to the recurrence given for (A001906) which is the sequence describing the sum. Though this is already a known fact, if we just bump the power up on $F_k$ to $F_k^3$, we still get a rather nice recurrence relation for the sum, in particular it is described by $0=(-N^4 + 7N^3- 9N^2-2N+4)x_n$. This integer sequence is recently added as number(A298591) in the OEIS. All powers of Fibonacci seem to follow this nice pattern that a linear recurrence where the terms do not depend on $n$ suffices, instead of in general for our technique, where the recurrence may need rational functions of $n$ showing up to describe the next term. These C-finite sequences are discussed in greater detail in \cite{Z2}. The techniques given in that paper can also be applied to some of the problems considered here.

Also of interest, suppose that $a_k$ is be the m-Fibonacci sequence, defined as $a_0=0$, $a_1=1$, and $a_{k+2} = ma_{k+1} + a_k$ for $k\ge 0$. Then, since the program was implemented in a symbolic way, there is no extra work to have that letting
\[
x_n = \sum_{k=0}^n a_k \binom{n}{k},
\]
we have 
\[
x_{n+2} = (2+m)x_{n+1} - x_{n}.
\]

This is also known, but is the main theorem of a twelve page paper by Falcon and Plaza \cite{FP}.

{\bf Application to multiple summations}
\\\\

Another promising application of our technique is to evaluating multiple sums over hypergeometric terms. A toy example of this would be computing

\[
\sum_{i=0}^n \sum_{k=0}^i \binom{i}{k} \binom{n}{i}.
\]

To find a recurrence for this sequence, pick out any of the factors which contain $k$, and run some automated process to evaluate single summation such as the Zeilberger Algorithm \cite{PWZ}. Often, this sum will not have a nice formula, so we are left with a possibly high order recurrence describing it. However, that is precisely what the techniques here are made to handle, so we can feed this partial evaluation into the procedure. Given enough computing this allows any number of summation signs to be dealt with. For each summation, we have the usual requirements of the original Sister Celine's method, namely that for each summation, the boundaries extend as far as the terms can be without becoming zero. In this particular case, evaluating the inner sum yields $0=(N - 2)x_n $, and substituting in that recurrence, we get that the whole sum satisfies $0=(N -3)z_n$. Which is to say, the sum evaluates to $3^n$. Though this has a nice combinatorial proof counting the number of assignments from $\{1,\ldots,n\}$ to $\{1,2,3\}$ by first picking the $k$ elements that map to either $1$ or $2$, and then, from those $k$ elements, picking the $i$ elements that map to $2$. That requires a moment of thought where such a simple recurrence for the computer only requires less than a second of `thought'. Alternatively, consider the harder problem, where we would want to compute

\[
\sum_{i=0}^n \sum_{k=0}^n \binom{i-k}{k}^2 \binom{n}{i}.
\]

It \textbf{may} be possible to evaluate this in a more human way, but for the computer it can easily determine that the solution is described by the recurrence

\[
0 = \left(-(n+9)N^5+ (7n+54)N^4-(17n+103)N^3+(21n+97)N^2-(15n+50)N+5n+5\right) x_n.
\]

A Maple package for multiple summations has already been described in \cite{AZ} and is available at:

{\tt http://sites.math.rutgers.edu/\~{}zeilberg/mamarim/mamarimhtml/multiZ.html}

However our package takes roughly the same time on the simple first example given, and is faster than their package on the second example. Their package, however, gives a `better' analysis of the summation, in that it does indefinite summation, and does not require that on the bounds of summation, the summand is zero. That is, theirs generalizes Zeilberger's algorithm, instead of Sister Celine's method.

{\bf Using this Maple package}
\\\\
Hopefully by this point, the usefulness of our package has been made clear. Though there is more detailed documentation in the maple package itself, here is a brief description of how they are used. The first step is to figure out the recurrence that is satisfied by $a_k$, called rec1. Then, call {\texttt findrec(I,J,timeout,rec1,F,d,n,N)} where both rec1 and the output are in shift operator notation, with N denoting the shift operator. This call will attempt to find the recurrence for the sum:

\[
x_n = \sum_{k=0}^n a_k^d H(n,k),
\]
where the recurrence is of order at most $I$, and degree at most $J$. {\texttt timeout} is the most time (in seconds) to wait on a particular attempt, if it exceeds that time, the procedure exits.

{\bf Acknowledgments}
\\\\
I would like to thank Doron Zeilberger for helping to guide me though this topic. I'd also like to thank Anthony Zaleski for his comments. The referee's detailed comments were invaluable to editing and clarifying this paper.

\pagebreak

\end{document}